\documentclass[12pt]{article}

\usepackage{amssymb,amsmath,amsthm, amsfonts}

\def\sqr#1#2{{\vcenter{\vbox{\hrule height.#2pt
              \hbox{\vrule width.#2pt height#1pt \kern#1pt \vrule width.#2pt}
              \hrule height.#2pt}}}}
\def\signed #1{{\unskip\nobreak\hfil\penalty50
              \hskip2em\hbox{}\nobreak\hfil#1
              \parfillskip=0pt \finalhyphendemerits=0 \par}}
\def\endpf{\signed {$\sqr69$}}
\def\3n{\negthinspace \negthinspace \negthinspace }
\def\2n{\negthinspace \negthinspace }
\def\1n{\negthinspace }

\def\={\buildrel \triangle \over =}

\def\a{\alpha}

\def\th{\theta}

\def\cL{{\cal L}}

\def\ms{\medskip}

\def\max{\mathop{\rm max}}
\def\min{\mathop{\rm min}}
\def\exp{\mathop{\rm exp}}
\def\sup{\mathop{\rm sup}}

\def\wt{\widetilde}
\def\cd{\cdot}

\def\inf{\hbox{\rm inf$\,$}}

\def\|{\Big |}
\def\({\Big (}
\def\){\Big )}
\def\[{\Big[}
\def\]{\Big]}
\def\be{\begin{equation}}
\def\bel{\begin{equation}\label}
\def\ee{\end{equation}}
\def\bt{\begin{theorem}}
\def\bcd{\begin{condition}}
\def\ecd{\end{condition}}
\def\et{\end{theorem}}
\def\bc{\begin{corollary}}
\def\ec{\end{corollary}}
\def\bde{\begin{definition}}
\def\ede{\end{definition}}
\def\bl{\begin{lemma}}
\def\el{\end{lemma}}
\def\bp{\begin{proposition}}
\def\ep{\end{proposition}}
\def\br{\begin{remark}}
\def\er{\end{remark}}
\def\ba{\begin{array}}
\def\ea{\end{array}}
\def\ed{\end{document}}

\def\square#1{\vbox{\hrule\hbox{\vrule height#1%
     \kern#1\vrule}\hrule}}
\def\rectangle#1#2{\vbox{\hrule\hbox{\vrule height#1%
     \kern#2\vrule}\hrule}}

\font\tenbb=msbm10 \font\sevenbb=msbm7 \font\fivebb=msbm5

\newfam\bbfam
\scriptscriptfont\bbfam=\fivebb \textfont\bbfam=\tenbb
\scriptfont\bbfam=\sevenbb

\newtheorem{lemma}{Lemma}[section]
\newtheorem{remark}{Remark}[section]

\newtheorem{theorem}{Theorem}[section]
\newtheorem{corollary}{Corollary}[section]

\newtheorem{definition}{Definition}[section]
\newtheorem{proposition}{Proposition}[section]
\newtheorem{condition}{Condition}[section]
\newtheorem{hypothesis}{Hypothesis}[section]

\makeatletter
   
   \@addtoreset{equation}{section}
\makeatother

\def\wt{\widetilde}

\begin{document}

\title{Multi-dimensional BSDE
with Oblique Reflection\\ and  Optimal Switching}

\author{Ying Hu\thanks{IRMAR, Universit\'e Rennes 1, Campus de Beaulieu,
35042 Rennes Cedex, France. Part of this work was completed when
this author was visiting Laboratory of Mathematics for Nonlinear
Sciences, Fudan University, whose hospitality is greatly
appreciated. {\small\it E-mail:} {\small\tt
Ying.Hu@univ-rennes1.fr}.\ms} \quad and \quad Shanjian
Tang\thanks{Institute of Mathematics and Department of Finance and
Control Sciences, School of Mathematical Sciences, Fudan
University, Shanghai 200433, China. This author  is supported in
part by NSFC Grant \#10325101, and the Chang Jiang Scholars
Programme. Part of this work was completed when this author was
visiting IRMAR, Universit\'e Rennes 1, whose hospitality is
greatly appreciated. {\small\it E-mail:} {\small\tt
sjtang@fudan.edu.cn}.\ms}}

\date{July 4, 2007}

\maketitle

\abstract{In this paper, we study a multi-dimensional backward
stochastic differential equation (BSDE) with oblique reflection,
which is a  BSDE reflected on the boundary of a special unbounded
convex domain along an oblique direction, and which arises
naturally in the study of optimal switching problem. The existence
of the adapted solution is obtained by the penalization method,
the monotone convergence, and the a priori estimations. The
uniqueness is obtained by a verification method (the first
component of any adapted solution is shown to be the vector value
of a switching problem for BSDEs). As applications, we apply the
above results to solve the optimal switching problem for
stochastic differential equations of functional type, and we give
also a probabilistic interpretation of the viscosity solution to a
system of variational inequalities.}
\medskip


\medskip

{\bf Abbreviated title. } Multi-dimensional  BSDEs with oblique
reflection


\section{Introduction}

In this paper, we are concerned with the following reflected
backward stochastic differential equation (RBSDE for short)  with
oblique reflection: for $i\in \Lambda:=\{1,\cdots,m\}$,

\be\label{RBSDEi}\left\{\ba{rcl}Y_i(t)&=&\displaystyle\xi_i+\int_t^T\psi(s,Y_i(s),Z_i(s),i)\,
ds-\int_t^TdK_i(s)-\int_t^TZ_i(s)\, dW(s),\ t\in [0,T],\\
Y_i(t)&\le&\displaystyle \min_{j\not=i}\{Y_j(t)+k(i,j)\}, \\
 &&\displaystyle
\int_0^T\left(Y_i(s)-\min_{j\not=i}\{Y_j(s)+k(i,j)\}\right)dK_i(s)=0.\ea\right.\ee
Here, $W$ is a standard Brownian motion on a complete probability
space $(\Omega, {\cal F},P)$ and $\xi$ is an $m$-dimensional
random variable measurable with respect to the past of $W$ up to
time $T$. $\xi$ is the terminal condition and $\psi$ the
coefficient (also called the generator). $k$ is a real function
defined on $\Lambda\times \Lambda$. The unknowns are the processes
$\{Y(t)\}_{t\in [0,T]}$, $\{Z(t)\}_{t\in [0,T]}$, and
$\{K(t)\}_{t\in [0,T]}$, which are required to be adapted with
respect to the natural completed filtration of the Brownian motion
$W$. Moreover, $K$ is an increasing process. The third relation in
(\ref{RBSDEi}) is called the minimal boundary condition.

\medskip

 RBSDE (\ref{RBSDEi}) evolves in the closure $\bar{Q}$ of domain $Q$:
$$Q:=\{(y_1,\cdots,y_m)^T\in \mathrm{R}^m: y_i< y_j+k(i,j) \ \mbox{\rm
for any } i,j\in \Lambda \mbox{ \rm such that }j\not=i\}, $$ which
is convex and unbounded. The boundary $\partial Q$ of domain $Q$
consists of the boundaries $\partial L_i^-, i\in\Lambda,$ with
$$L_{i}^-:=\{(y_1,\cdots,y_m)^T\in \mathrm{R}^m:y_i<y_j+k(i,j),
 \mbox{ for any }j\in\Lambda \mbox{ such that } j\not=i\}, i\in\Lambda.$$
That is,
$$\partial Q=\mathop{\cup}_{i=1}^m\partial L_i^-. $$
In the interior of $\bar{Q}$, each equation in (\ref{RBSDEi}) is
independent of others. On the boundary, say $\partial L_{i}^-$,
the $i$-th equation  is switched to another one, and the solution
is reflected along the oblique direction $-e_i$ (which is the
opposite direction of the $i$-th coordinate axis).

RBSDE was first studied by El Karoui et al. \cite{ElKPPQ} for the
one-dimensional case. Multidimensional RBSDE was studied by
Gegout-Petit and Pardoux \cite{Gegout-PetitPardoux}, but their
BSDE is reflected on the boundary of a convex domain along the
inward normal direction, and their method depends heavily on the
properties of this inward normal reflection (see (1)-(3) in
\cite{Gegout-PetitPardoux}). We note that in a very special case
(e.g., $\psi$ is independent of $z$), Ramasubramanian
\cite{Ramasubramanian} studied a BSDE in an orthant with oblique
reflection. Note also that there are some papers dealing with SDEs
with oblique reflection (see, e.g. \cite{LionsSznitman} and
\cite{DupuisIshii}).

An incomplete and less general form of RBSDE (\ref{RBSDEi}) (where
the minimal condition of (\ref{RBSDEi}) is missing and the
generator $\psi$ does not depend on $(y,z)$) is suggested by
\cite{CarLud}. But they did not discuss the existence and
uniqueness of solution, which is considered to be difficult. See
Remark 3.1 in \cite{CarLud}.

\medskip

Besides the theoretic interest, RBSDE (\ref{RBSDEi}) arises
naturally from the following optimal switching problem.

Consider the switched equation
\begin{equation}\label{sdef1}
X^{a(\cdot)}(t)=x_0+\int_0^t
\sigma(s,X^{a(\cdot)})[dW(s)+b(s,X^{a(\cdot)},a(s))ds], \ t\in
[0,T]
\end{equation}
and the cost functional
\begin{equation}\label{cost1}
J(a(\cdot))=E[\int_0^T
l(s,X^{a(\cdot)},a(s))ds]+E[\sum_{i=1}^\infty
k(\alpha_{i-1},\alpha_i)].
\end{equation}
The optimal switching problem is to minimize the cost
$J(a(\cdot))$ with respect to $a(\cdot)$, subject to the state
equation (\ref{sdef1}).

 In the above, $x_0$ is a fixed point in $R^d$. $\sigma$, $b$ and $l$ are defined on
$[0,T]\times C([0,T];R^d)$, $[0,T]\times C([0,T];R^d)\times
\Lambda$ and $[0,T]\times C([0,T];R^d)\times \Lambda$,
respectively, with appropriate dimensions.
$$a(\cdot)=\alpha_{0}\chi_{\{\th_0\}}(\cdot)+\sum_{i=1}^{\infty} \alpha_{i-1}\chi_{(\th_{i-1},\th_i]}(\cdot)$$
is called an admissible switching strategy if for any $i$,
$\theta_i$ is a stopping time, and $\alpha_i$ is an ${\cal
F}_{\th_i}$-measurable random variable with values in $\Lambda$.
Here, $\chi$ is the indicator function. $k$ is called the
switching cost.

\medskip

Optimal switching is a special case of impulse control. The
classical method of quasi-variational inequalities to solve
impulse control problems driven by Markov processes is referred to
the book of Bensoussan and Lions \cite{BL}. See  \cite{TY} and the
references therein for the theory of variational inequalities and
the dynamic programming for optimal stochastic switching. But
these works are restricted within the Markovian case. Recently,
using the method of  Snell envelope  (see, e.g. El Karoui
\cite{ElK}) combined with the theory of  scalar valued RBSDEs,
Hamadene and Jeanblanc \cite{HaJean} studied the switching problem
with two modes (i.e., $m=2$) in the non-Markovian context.
Djehiche, Hamadene and Popier \cite{DHP} generalized their result
to the above switching problem with multi modes. We note that in
both \cite{HaJean} and \cite{DHP}, the drift term $b$ does not
depend on the switching strategy $a$.

\medskip

The main contribution of this paper is to establish the existence
and uniqueness of solution for RBSDE (\ref{RBSDEi}). We prove the
existence by the penalization method, the monotone convergence,
and the a priori estimation whose proof is rather technical. The
proof of uniqueness is quite different: the classical method to
estimate the difference of two solutions appears difficult to be
applied to our present case of the  oblique reflection. We obtain
the uniqueness by a verification method: first we introduce an
optimal switching problem for BSDEs, then we prove that the first
component $Y$ of any adapted solution $(Y,Z,K)$ of RBSDE
(\ref{RBSDEi}) is the (vector) value for the switching problem. As
applications, we solve the optimal switching problem (\ref{sdef1})
and (\ref{cost1}), and we establish the Feynman-Kac formula for
the viscosity solution to a system of variational inequalities.

\medskip

The paper is organized as follows: in Section 2, we prove the
existence of solution, whereas Section 3 is devoted to the
uniqueness. We solve the optimal switching problem in Section 4.
Finally, in Section 5, we give a probabilistic interpretation of
the viscosity solution  to a system of variational inequalities.

\section{Existence}

\subsection{Notations}

Let us fix a nonnegative real number $T>0$. First of all,
$W=\{W_t\}_{t\ge 0}$ is a standard Brownian motion with values in
$R^d$ defined on some complete probability space $(\Omega,{\cal
F},P)$. $\{{\cal F}_t,t\ge 0\}$ is the natural filtration of the
Brownian motion $W$ augmented by the $P$-null sets of ${\cal F}$.
All the measurability notions will refer to this filtration. In
particular, the sigma-field of predictable subsets of $[0,T]\times
\Omega$ is denoted by ${\cal P}$.

Let us consider now the RBSDE (\ref{RBSDEi}). The generator $\psi$
is a random function $\psi:[0,T]\times \Omega\times R\times
R^d\times \Lambda \rightarrow R$ whose component $\psi(\cdot,i)$
is measurable with respect to ${\cal P}\otimes {\cal B}(R)\otimes
{\cal B}(R^d)$ and the terminal condition $\xi$ is simply a
$R^m$-valued ${\cal F}_T$-measurable random variable. $k$ is
defined on $\Lambda\times \Lambda$ and scalar valued.

\medskip

By a solution to  RBSDE (\ref{RBSDEi}) we mean a triple
$(Y,Z,K)=\{Y(t),Z(t),K(t)\}_{t\in [0,T]}$ of predictable processes
with values in $R^m\times R^{m\times d}\times R^m$ such that
$P$-a.s., $t\rightarrow Y(t)$ and $t\rightarrow K(t)$ are
continuous, $t\rightarrow Z(t)$ belongs to $L^2(0,T)$,
$t\rightarrow \psi(t,Y_i(t),Z_i(t),i)$ belongs to $L^1(0,T)$ and
$P$-a.s.,  RBSDE (\ref{RBSDEi}) holds.

$S^2(R^m)$ or simply $S^2$ denotes the set of $R^m$-valued,
adapted and c\`adl\`ag processes $\{Y(t)\}_{t\in [0,T]}$ such that
$$||Y||_{S^2}:=E[\sup_{t\in [0,T]}|Y(t)|^2]^{1/2}<+\infty.$$
$(S^2,||\cdot||_{S^2})$ is a Banach space.

$M^2(R^{m\times d})$ or simply $M^2$ denotes the set of
(equivalent classes of) predictable processes $\{Z(t)\}_{t\in
[0,T]}$ with values in $R^{m\times d}$ such that
$$||Z||_{M^2}:=E[\int_0^T |Z(s)|^2ds]^{1/2}<+\infty.$$
$M^2$ is a Banach space endowed with this norm.

\begin{eqnarray*}
N^2(R^m):&=&\{K=(K_1,\cdots, K_m)^T\in S^2 :\mbox{ for any } i\in
\Lambda, K_i(0)=0, \\
& &\mbox{ and } t\rightarrow K_i(t) \mbox{ is increasing }\},
\end{eqnarray*}
where $^T$ means transpose. $(N^2,||\cdot||_{S^2})$ is a Banach
space.

\subsection{Existence}
In this subsection, we prove the existence result for  RBSDE
(\ref{RBSDEi}). We assume the following Lipschiz condition on the
generator.
\begin{hypothesis}\label{psi}
(i)
$\psi(\cdot,0,0):=(\psi(\cdot,0,0,1),\cdots,\psi(\cdot,0,0,m))^T$
belongs to $M^2$.

(ii) There exists a constant $C> 0$, such that, $P$-a.s. for each
$(t,y,y',z,z',i)\in [0,T]\times R\times R\times R^{d}\times
R^{d}\times\Lambda$,
$$|\psi(t,y,z,i)-\psi(t,y',z',i)|\le C (|y-y'|+|z-z'|).$$
\end{hypothesis}

We make the following assumption on $k$ which is standard in the
literature of optimal switching.

\begin{hypothesis}\label{k}
(i) For any $(i,j)\in \Lambda\times \Lambda$, $k(i,j)\ge 0$.

(ii) For any $(i,j,l)\in \Lambda\times \Lambda\times \Lambda$,
$$k(i,j)+k(j,l)\ge k(i,l).$$
\end{hypothesis}

We are now in position to state the existence result.
\begin{theorem}\label{existence}
Let the Hypotheses \ref{psi} and \ref{k} hold. Assume that $\xi\in
L^2(\Omega,{\cal F}_T,P;R^m)$ takes values in $\bar{Q}$. Then
RBSDE (\ref{RBSDEi}) has a solution $(Y,Z,K)$ in $S^2\times
M^2\times N^2$.
\end{theorem}

We first sketch our proof.

{\bf Sketch of the Proof:} The proof is divided to four steps. In
Step 1, we introduce the penalized BSDEs whose existence and
uniqueness follows from the classical result. In Step 2, we state
some (uniform) a priori estimates for the solutions of penalized
BSDEs, whose proof will be given in the next subsection. In Step
3, we prove the (monotone) convergence of these solutions.
Finally, in Step 4, we check out the minimal boundary condition.

{\bf Proof of Theorem \ref{existence} :} Step 1. The penalized
BSDEs.

\medskip

For any nonnegative integer $n$, let us introduce  the following
penalized BSDE:
\begin{eqnarray}\label{penalized}\begin{array}{rcl}
Y_i^n(t)&=&\displaystyle \xi_i
+\int_t^T\psi(s,Y_i^n(s),Z_i^n(s),i)\,
ds \\[0.3cm]
&&\displaystyle
-n\sum_{l=1}^m\int_t^T(Y_i^n(s)-Y_l^n(s)-k(i,l))^+\,
ds\\[0.5cm]
&&\displaystyle-\int_t^T Z_i^n(s)\, dW(s),\ t\in [0,T],\ i\in
\Lambda.
\end{array}\label{PBSDE}\end{eqnarray}
Note that when $l=i$, we have, in view of Hypothesis \ref{k} (i),
\be(Y_i^n(s)-Y_l^n(s)-k(i,l))^+=0. \ee

From the classical result of Pardoux and Peng \cite{PP}, for any
$n$, BSDE (\ref{penalized}) has a unique solution $(Y^n,Z^n)$ in
the space $S^2\times M^2$.

\medskip

Step 2. A priori estimates.

The following lemma will play a crucial rule in the proof of
Theorem \ref{existence}.
\begin{lemma}\label{apriori}
Let the Hypotheses \ref{psi} and \ref{k} hold. Let us also assume
that $\xi\in L^2(\Omega,{\cal F}_T,P;R^m)$ takes values in
$\bar{Q}$. Then there exists a constant $C>0$ (independent of
$n$), such that
\begin{equation}\label{estimateyz}
||Y^n||_{S^2}+||Z^n||_{M^2}\le C
\end{equation}
and \begin{equation}\label{estimateK}
n^2E\int_0^T\left((Y_i^n(s)-Y_j^n(s)-k(i,j))^+\right)^2\, ds\le C.
\end{equation}
\end{lemma}

However, the proof of this lemma is quite lengthy and delicate. We
relegate it to the next subsection.

\medskip

Step 3. Convergence of  solutions $\{Y^n,Z^n\}$ of penalized
BSDEs.

First, for each $n$, we introduce a function $\psi^n$ as follows:
$$\psi^n(t,y,z,i):=\psi(t,y_i,z_i,i)-n\sum_{l=1}^m
(y_i-y_l-k(i,l))^+, \quad (t,y,z,i)\in [0,T]\times R^m\times
R^{m\times d}\times\Lambda.$$

Since $\psi^n(\cdot,z,i)$ depends only on $z_i$ and
$$\psi^n(t,y+y',z_i,i)\ge \psi^n(t,y',z_i,i) $$
for any $y\in R^m \mbox{ such that } y\ge 0 \mbox{ and the
}i-\mbox{th component }  y_i=0,$ it is easy to check
$$-4\langle y^-,\psi^n(t,y^+ +y',z)-\psi^{n+1}(t,y',z')\rangle \le 2\sum_{i=1}^m
\chi_{y_i<0}|z_i-z_i'|^2+C|y'|^2, \quad P-a.s.,$$ for a constant
$C>0$. We note that this is the inequality (5) in \cite{HuPeng}.
Applying the comparison theorem for multi-dimensional BSDEs (see
Hu and Peng \cite [Theorem 2.1] {HuPeng}), we deduce that for any
nonnegative integer $n$,
 \be  Y_i^n(t)\ge Y_i^{n+1}(t), \forall i\in \Lambda, t\in [0,T]. \ee

For a.e. $t$ and $P$-a.s. $\omega$, $\{Y^n(t,\omega)\}_n$ admits a
limit, denoted by $Y(t,\omega)$. Moreover from the a priori
estimate (\ref{estimateyz}) and Fatou's lemma, we have \be
\sup_{t\in [0,T]} E|Y(t)|^2\le C.\ee

In view of the fact  that $Y_i(t)\le Y_i^n(t)\le Y_i^0(t)$ with
$i\in\Lambda$,  applying Lebesgue's dominated convergence theorem,
we have
 \be\label{ynm}
\lim_{n\rightarrow\infty} E\int_0^T|Y^n(s)-Y(s)|^2\, ds= 0. \ee

Now we prove that $\left\{(Y^n,Z^n)\right\}_n$ is a Cauchy
sequence in the space $S^2\times M^2$. For this purpose, we apply
It\^o's formula to $|Y^n_i(t)-Y^p_i(t)|^2$ to obtain

\begin{eqnarray}\label{itoynm}
& &|Y_i^n(t)-Y_i^p(t)|^2+\int_t^T|Z_i^n(s)-Z_i^p(s)|^2\,
ds\nonumber\\ &=&2\int_t^T(Y_i^n(s)-Y_i^p(s))(\psi(s,Y_i^n(s),
Z_i^n(s),i)-\psi(s,Y_i^p(s),Z_i^p(s),i))\,
ds\nonumber \\
& & -2\int_t^T(Y_i^n(s)-Y_i^p(s))n(Y_i^n(s)-Y_j^n(s)-k(i,j))^+\, ds\nonumber\\
& &+2\int_t^T(Y_i^n(s)-Y_i^p(s))p(Y_i^p(s)-Y_j^p(s)-k(i,j))^+\, ds\nonumber\\
& & -2\int_t^T(Y_i^n(s)-Y_i^p(s))(Z_i^n(s)-Z_i^p(s))\, dW(s), i\in
\Lambda.
\end{eqnarray}
Putting $t=0$ and taking expectation in the last equality, we get
for $i\in\Lambda$,

\begin{eqnarray*}
 & &E|Y_i^n(0)-Y_i^p(0)|^2+E\int_0^T|Z_i^n(s)-Z_i^p(s)|^2 ds\\
 &=& 2E \int_0^T(Y_i^n(s)-Y_i^p(s))
 (\psi(s,Y_i^n(s), Z_i^n(s),i)-\psi(s,Y_i^p(s),Z_i^p(s),i)) ds\\
 & &-2E\int_0^T(Y_i^n(s)-Y_i^p(s))n(Y_i^n(s)-Y_j^n(s)-k(i,j))^+
 ds\\
 & &+2E\int_0^T(Y_i^n(s)-Y_i^p(s))p(Y_i^p(s)-Y_j^p(s)-k(i,j))^+  ds\\
 &\le & C E\int_0^T(Y_i^n(s)-Y_i^p(s))^2 ds
 +{\frac{1}{2}}E\int_0^T|Z_i^n(s)-Z_i^p(s)|^2 ds\\
 & &+\left(E\int_0^T|Y_i^n(s)-Y_i^p(s)|^2\, ds\right)^{\frac{1}{2}} \left(E\int_0^Tn^2\left((Y_i^n(s)-Y_j^n(s)-k(i,j))^+\right)^2
 ds\right)^{\frac{1}{2}}\\
 & &+\left(E\int_0^T|Y_i^n(s)-Y_i^p(s)|^2\, ds\right)^{\frac{1}{2}}\left(E\int_0^Tp^2\left((Y_i^p(s)-Y_j^p(s)-k(i,j))^+\right)^2
ds\right)^{\frac{1}{2}}.
\end{eqnarray*}
From (\ref{estimateK}) and (\ref{ynm}), we have
 \be\label{znm} \lim_{n,p\rightarrow \infty}E\int_0^T|Z_i^n(s)-Z_i^p(s)|^2\, ds= 0,i\in\Lambda. \ee

Now, we define the increasing process $K_i^n$ as follows:
\begin{equation}\label{Kn}
K_i^n(t):=n\int_0^t\sum_{l=1}^m(Y_i^n(s)-Y_l^n(s)-k(i,l))^+\, ds,
\quad t\in [0,T],i\in \Lambda.
\end{equation}
From the penalized BSDE (\ref{penalized}), we have
\begin{equation}\label{Kin}
K_i^n(t)=Y_i^n(t)-Y_i^n(0)+\int_0^t
\psi(s,Y_i^n(s),Z^n_i(s),i)ds-\int_0^t Z^n_i(s)dW(s),i\in\Lambda.
\end{equation}

We denote by $Z$  the limit of $Z^n$ in $M^2$. Set
\begin{equation}\label{eqK}
K_i(t):=Y_i(t)-Y_i(0)+\int_0^t \psi(s,Y_i(s),Z_i(s),i)ds-\int_0^t
Z_i(s)dW(s),i\in\Lambda.
\end{equation}
We have
$$\lim_{n\to \infty} ||Z^n-Z||_{M^2}=0.$$

Going back again to (\ref{itoynm}), we deduce that for
$i\in\Lambda$,

\begin{eqnarray}\label{supynm}
& &E\sup_{0\le t\le T}|Y_i^n(t)-Y_i^p(t)|^2 \nonumber\\
&\le&2E\int_0^T|Y_i^n(s)-Y_i^p(s)||\psi(s,Y_i^n(s),
Z_i^n(s),i)-\psi(s,Y_i^p(s),Z_i^p(s),i)|\, ds\nonumber\\
& & +2E\int_0^T|Y_i^n(s)-Y_i^p(s)|n(Y_i^n(s)-Y_j^n(s)-k(i,j))^+\,
ds\nonumber\\
& &
 +2E\int_0^T|Y_i^n(s)-Y_i^p(s)|p(Y_i^p(s)-Y_j^p(s)-k(i,j))^+\, ds\nonumber\\
& &+2E\sup_{0\le t\le
T}|\int_t^T(Y_i^n(s)-Y_i^p(s))(Z_i^n(s)-Z_i^p(s))\, dW(s)|.
\end{eqnarray}

 The last term of the right hand side of the last inequality is less than or equal to the following
quantity: \begin{eqnarray}\label{lastterm} & &
CE\left(\int_0^T|(Y_i^n(s)-Y_i^p(s))(Z_i^n(s)-Z_i^p(s))|^2\,
ds\right)^{\frac{1}{2}}\nonumber \\
&\le& CE\left(\sup_{0\le t\le
T}|Y_i^n(t)-Y_i^p(t)|\int_0^T|Z_i^n(s)-Z_i^p(s)|^2\,
ds\right)^{\frac{1}{2}}\nonumber \\
&\le &{\frac{1}{2}}E\sup_{0\le t\le
T}|Y_i^n(t)-Y_i^p(t)|^2+CE\int_0^T|Z_i^n(s)-Z_i^p(s)|^2\, ds.
\end{eqnarray}

Combining (\ref{supynm}) and (\ref{lastterm}) and taking into
consideration (\ref{ynm}) and (\ref{znm}), we deduce  that
$\{Y^n\}_n$ is a Cauchy sequence in $S^2$, which means in
particular that
$$\lim_{n\rightarrow\infty}||Y^n-Y||_{S^2}= 0.$$
Consequently, $Y$ is a continuous process.

From (\ref{Kin}), (\ref{eqK}),  and  the following fact that
$$||Y^n-Y||_{S^2}+||Z^n-Z||_{M^2}\rightarrow 0,\quad\mbox{as } n\to \infty,$$
 we deduce immediately that
$$\lim_{n\to\infty}||K^n-K||_{S^2}=0.$$
Hence,  $K \in N^2$, and $(Y,Z,K)$ satisfies the first relation in
RBSDE (\ref{RBSDEi}).

Finally, from the a priori estimate (\ref{estimateK}), we have
$$E\int_0^T\left((Y_i^n(s)-Y_j^n(s)-k(i,j))^+\right)^2\, ds\le
\frac{C}{n^2}, \ i,j\in\Lambda.$$ Letting $n\rightarrow \infty$,
we deduce
$$E\int_0^T\left((Y_i(s)-Y_j(s)-k(i,j))^+\right)^2\, ds=0, \ i,j\in\Lambda.$$ Hence,
\be\label{inQ} Y_i(s)\le Y_j(s)+k(i,j), \quad s\in [0,T],
i,j\in\Lambda,\ee which are equivalent to the following:
$$P-a.s. \quad Y(s)\in \bar{Q}, \quad \forall s\in [0,T].$$

Step 4. The minimal boundary condition.

 Let us first state the following lemma whose proof is at the end of this subsection.
\bl\label{ynkn} Let the Hypotheses \ref{psi} and \ref{k} hold. Let
us also assume that $\xi\in L^2(\Omega,{\cal F}_T,P;R^m)$ takes
values in $\bar{Q}$. We have, for any integer $n$,
 \be\label{Kn0} \int_0^T\left(Y_i^n(s)-\min_{j\not=i}[Y_j^n(s)+k(i,j)]\right)^-\, dK^n_i(s)=0,
 \quad i\in \Lambda.\ee
 \el

Now, we can take the limit in (\ref{Kn0}) by letting $n$ tend to
$+\infty$ and applying Lemma 5.8 in \cite{Gegout-PetitPardoux} to
get the following \be
\int_0^T\left(Y_i(s)-\min_{j\not=i}\left[Y_j(s)+k(i,j)\right]\right)^-\,
dK_i(s)=0, \quad i\in\Lambda,\ee which, together with (\ref{inQ}),
yields the minimal boundary conditions.

The proof of Theorem \ref{existence} is now complete.
\endpf

{\bf Proof of Lemma \ref{ynkn}}  For $i\in\Lambda$, the left hand
side of the last equality is equal to the following sum \be n
\sum_{l=1}^m\int_0^T\min_{j\not=i}\left\{\left(Y_i^n(s)-Y_j^n(s)-k(i,j)\right)^-
(Y_i^n(s)-Y_l^n(s)-k(i,l))^+\right\}\, ds.  \ee

We claim that the integrand of the $l$-th integral is equal to
zero for $l\in\Lambda$.

In fact, it is immediate for the case of $l=i$. For the case of
$l\not=i$, the integrand is the minimum of the following $m-1$
nonnegative quantities: \be\left(Y_i^n(s)-Y_j^n(s)-k(i,j)\right)^-
\left(Y_i^n(s)-Y_l^n(s)-k(i,l)\right)^+, \quad j\in\Lambda \mbox{
\rm and } j\not=i,\ee whose $l$-th term is zero due to the fact
that \be \left(Y_i^n(s)-Y_l^n(s)-k(i,l)\right)^-
\left(Y_i^n(s)-Y_l^n(s)-k(i,l)\right)^+=0, \ee and therefore, it
is zero. The desired result then follows.
\endpf

\subsection{Proof of Lemma \ref{apriori}}

In this subsection, we prove Lemma \ref{apriori}.

For $i,j\in\Lambda$, applying Tanaka's formula (see, e.g.
\cite{RY}) to $(Y_i^n(t)-Y_j^n(t)-k(i,j))^+$, we have
 \be\ba{rcl}
&&\displaystyle
(Y_i^n(t)-Y_j^n(t)-k(i,j))^++n\sum_{l=1}^m\int_t^T\chi_{\cL_{ij,n}^+}(s)(Y_i^n(s)-Y_l^n(s)-k(i,l))^+\,
ds\\
&&\displaystyle-n\sum_{l=1}^m\int_t^T\chi_{\cL_{ij,n}^+}(s)(Y_j^n(s)-Y_l^n(s)-k(j,l))^+\,
ds+{\frac{1}{2}}\int_t^T\, dL_{ij}^n(s)\\
&=&\displaystyle \int_t^T\chi_{\cL_{ij,n}^+}(s) (\psi(s,Y_i^n(s),
Z_i^n(s),i)-\psi(s,Y_j^n(s), Z_j^n(s), j))\, ds
\\
&&\displaystyle
-\int_t^T\chi_{\cL_{ij,n}^+}(s)(Z_i^n(s)-Z_j^n(s))\, dW(s), \ea
\ee where for $i,j\in\Lambda$,\be \ba{c}
\cL_{ij,n}^+:=\{(s,\omega):Y_i^n(s)>Y_j^n(s)+k(i,j)\}, \ea \ee and
$L_{ij}^n$ is the local time of the process $Y_i^n-Y_j^n-k(i,j)$
at $0$.

Applying It\^o's formula to
$\left((Y_i^n(t)-Y_j^n(t)-k(i,j))^+\right)^2$ and taking into
consideration \be \int_t^T(Y_i^n(s)-Y_j^n(s)-k(i,j))^+\,
dL_{ij}^n(s)=0,\forall t\in [0,T],\ee we have
  \be\label{Tanaka}\ba{rcl}
&&\displaystyle\left((Y_i^n(t)-Y_j^n(t)-k(i,j))^+\right)^2\\
&& \displaystyle
+2n\int_t^T\left((Y_i^n(s)-Y_j^n(s)-k(i,j))^+\right)^2\,
ds\\[0.5cm]
&&\displaystyle+\int_t^T\chi_{\cL_{ij,n}^+}(s)|Z_i^n(s)-Z_j^n(s)|^2\, ds\\[0.5cm]
&=&\displaystyle 2\int_t^T
(Y_i^n(s)-Y_j^n(s)-k(i,j))^+\left[\psi(s,Y_i^n(s),
Z_i^n(s),i)-\psi(s,Y_j^n(s), Z_j^n(s),j)\right]\,
ds\\[0.5cm]
&&\displaystyle-2\int_t^T(Y_i^n(s)-Y_j^n(s)-k(i,j))^+(Z_i^n(s)-Z_j^n(s))\,
dW(s)\\
&&\displaystyle
+2n\int_t^T(Y_i^n(s)-Y_j^n(s)-k(i,j))^+(Y_j^n(s)-Y_i^n(s)-k(j,i))^+\,
ds\\[0.5cm]
&&\displaystyle +2n\sum_{l\not=
i,l\not=j}\int_t^T(Y_i^n(s)-Y_j^n(s)-k(i,j))^+\\
&&\displaystyle \times [(Y_j^n(s)-Y_l^n(s)-k(j,l))^+
-(Y_i^n(s)-Y_l^n(s)-k(i,l))^+]\, ds. \ea \ee

We claim that the integrands of the integrals in the last two
terms of (\ref{Tanaka}) are all less than or equal to zero. In
fact, since
$$\{(y_1,\cdots,y_m)^T\in R^m: y_i-y_j-k(i,j)>0,y_j-y_i-k(j,i)>0\}=\emptyset$$
due to the fact that $$k(i,j)+k(j,i)\ge 0,$$ we have
 \be
(Y_i^n(s)-Y_j^n(s)-k(i,j))^+(Y_j^n(s)-Y_i^n(s)-k(j,i))^+ =0, \quad
i,j\in\Lambda.\ee

On the other hand, for $l,i,j\in\Lambda$,  taking into
consideration both Hypothesis \ref{k} (ii), i.e.,
$$k(i,j)+k(j,l)\ge k(i,l),$$
and the elementary inequality that $x_1^+-x_2^+\le (x_1-x_2)^+$
for any two real numbers $x_1$ and $x_2$, we have \be \ba{rcl}
&&(Y_i^n(s)-Y_j^n(s)-k(i,j))^+[(Y_j^n(s)-Y_l^n(s)-k(j,l))^+
-(Y_i^n(s)-Y_l^n(s)-k(i,l))^+]\\
 &\le&
(Y_i^n(s)-Y_j^n(s)-k(i,j))^+(Y_j^n(s)-Y_i^n(s)-k(j,l)+k(i,l))^+\\
 &\le
&(Y_i^n(s)-Y_j^n(s)-k(i,l)+k(j,l))^+(Y_j^n(s)-Y_i^n(s)-k(j,l)+k(i,l))^+.\ea\ee
The last term of the last inequality is zero, since $$
\{(y_1,\cdots,y_m)^T\in R^m:y_i-y_j-k(i,l)+k(j,l)>0,
y_j-y_i-k(j,l)+k(i,l)>0\}=\emptyset.$$

Concluding the above, we have \be\label{TanakaE}\ba{rcl}
&&\displaystyle E\left((Y_i^n(t)-Y_j^n(t)-k(i,j))^+\right)^2\\[0.5cm]
&& \displaystyle
+2nE\int_t^T\left((Y_i^n(s)-Y_j^n(s)-k(i,j))^+\right)^2\,
ds\\[0.5cm]
&&\displaystyle+E\int_t^T\chi_{\cL_{ij,n}^+}(s)|Z_i^n(s)-Z_j^n(s)|^2\, ds\\[0.5cm]
&\le& \displaystyle 2E\int_t^T (Y_i^n(s)-Y_j^n(s)-k(i,j))^+
|\psi(s,Y_i^n(s), Z_i^n(s),i)-\psi(s,Y_j^n(s), Z_j^n(s),j)|\,
ds.\ea \ee

In view of Hypothesis \ref{psi} on the  function $\psi$,  we
have \be\ba{rcl} &&|\psi(s,Y_i^n(s), Z_i^n(s),i)-\psi(s,Y_j^n(s), Z_j^n(s),j)|\\
 &\le& |\psi(s,Y_i^n(s), Z_i^n(s),i)-\psi(s,Y_i^n(s),
 Z_i^n(s),j)|\\
 & &
+|\psi(s,Y_i^n(s), Z_i^n(s),j)-\psi(s,Y_j^n(s), Z_j^n(s),j)|\\
&\le&
C(|\psi(s,0,0)|+|Y_i^n(s)|+|Z_i^n(s)|+|Y_i^n(s)-Y_j^n(s)|+|Z_i^n(s)-Z_j^n(s)|)\\
 &\le&
C(1+|\psi(s,0,0)|+|Y_i^n(s)|+|Z_i^n(s)|\\
& &+|Y_i^n(s)-Y_j^n(s)-k(i,j)|+|Z_i^n(s)-Z_j^n(s)|).\ea\ee

Consequently, we have
 \be\ba{rcl}
&&\displaystyle E\left((Y_i^n(t)-Y_j^n(t)-k(i,j))^+\right)^2\\[0.5cm]
&&
\displaystyle+2nE\int_t^T\left((Y_i^n(s)-Y_j^n(s)-k(i,j))^+\right)^2\,
ds\\[0.5cm]
&&\displaystyle+E\int_t^T\chi_{\cL_{ij,n}^+}(s)|Z_i^n(s)-Z_j^n(s)|^2\, ds\\[0.5cm]
&\le& \displaystyle C E\int_t^T
|(Y_i^n(s)-Y_j^n(s)-k(i,j))^+|^2\, ds\\[0.5cm]
&&\displaystyle+{\frac{1}{2}}E\int_t^T\chi_{\cL_{ij,n}^+}(s)(1+|\psi(s,0,0)|^2+|Y_i^n(s)|^2+|Z_i^n(s)|^2\\
& &+|(Y_i^n(s)-Y_j^n(s)-k(i,j))^+|^2+|Z_i^n(s)-Z_j^n(s)|^2)\,
ds.\ea \ee

Applying Gronwall's inequality, we deduce easily that
\be\ba{c}\displaystyle
E\left((Y_i^n(t)-Y_j^n(t)-k(i,j))^+\right)^2\le
C\left(1+E\int_0^T\chi_{\cL_{ij,n}^+}(s)\left(|Y_i^n(s)|^2+|Z_i^n(s)|^2\right)\, ds \right),\\[0.5cm]
\displaystyle
nE\int_0^T\left((Y_i^n(s)-Y_j^n(s)-k(i,j))^+\right)^2\, ds
+E\int_0^T\chi_{\cL_{ij,n}^+}(s)|Z_i^n(s)-Z_j^n(s)|^2\, ds\\[0.5cm]
\displaystyle \le
C\left(1+E\int_0^T\chi_{\cL_{ij,n}^+}(s)[|Y_i^n(s)|^2+|Z_i^n(s)|^2]\,
ds \right). \ea\ee

Going back to (\ref{Tanaka}) and applying Burkholder-Davis-Gundy's
inequality, we obtain \be\label{sup} \displaystyle E[\sup_{0\le
t\le T}\left((Y_i^n(t)-Y_j^n(t)-k(i,j))^+\right)^2]\le
C\left(1+E\int_0^T\chi_{\cL_{ij,n}^+}(s)[|Y_i^n(s)|^2+|Z_i^n(s)|^2]\,
ds \right).\ee

On the other hand, from (\ref{TanakaE}), we deduce that,
\begin{eqnarray*} & & E\left((Y_i^n(t)-Y_j^n(t)-k(i,j))^+\right)^2
+2nE\int_t^T\left((Y_i^n(s)-Y_j^n(s)-k(i,j))^+\right)^2\,
ds\\
& &+E\int_t^T\chi_{\cL_{ij,n}^+}(s)|Z_i^n(s)-Z_j^n(s)|^2\, ds\\
&\le&  (n+C)E\int_t^T
\left((Y_i^n(s)-Y_j^n(s)-k(i,j))^+\right)^2\, ds\\
&
&+{\frac{C}{n}}E\int_t^T\chi_{\cL_{ij,n}^+}(s)(1+|\psi(s,0,0)|^2+|Y_i^n(s)|^2+|Z_i^n(s))|^2+|Z_i^n(s)-Z_j^n(s)|^2)\,
ds. \end{eqnarray*} This shows that, for sufficiently large $n$,
 \be\label{nn}\displaystyle
n^2E\int_0^T\left((Y_i^n(s)-Y_j^n(s)-k(i,j))^+\right)^2\, ds
\le C\left(1+E\int_0^T[|Y_i^n(s)|^2+|Z_i^n(s)|^2]\, ds \right).\\
\ee

Finally, applying It\^o's formula to $|Y_i^n(t)|^2$, we obtain:
\begin{eqnarray}\label{finalito}
& &|Y_i^n(t)|^2+\int_t^T|Z_i^n(s)|^2\, ds\nonumber \\
&=& |\xi_i|^2+2\int_t^T Y_i^n(s)\cdot
\left[\psi(s,Y_i^n(s),Z_i^n(s),i)-\sum_{l=1}^m
n(Y_i^n(s)-Y_l^n(s)-k(i,l))^+\right]\,
 ds\nonumber\\
 & &-2\int_t^T Z_i(s)dW(s).
\end{eqnarray}

By taking expectation and using the elementary inequality:
$$2ab\le \frac{1}{\epsilon}a^2+\epsilon b^2,\quad \forall \epsilon>0,$$ we deduce that, for any
$\epsilon>0$,
 \begin{eqnarray*} & &E|Y_i^n(t)|^2+E\int_t^T|Z_i^n(s)|^2\, ds\\
& & \le E|\xi_i|^2+2E\int_t^T|Y_i^n(s)|\cdot
 \bigr[|\psi(s,Y_i^n(s),Z_i^n(s),i)|\\
& &+\sum_{l=1}^m n(Y_i^n(s)-Y_l^n(s)-k(i,l))^+\bigr]\,
 ds\\
&\le& C+2E\int_t^T|Y_i^n(s)|\cdot
\Bigr[C(|\psi(s,0,0)|+|Y_i^n(s)|+|Z_i^n(s)|)\\
& &+\sum_{l=1}^m n(Y_i^n(s)-Y_l^n(s)-k(i,l))^+\Bigr]\,
 ds\\
& \le &C+C_\epsilon E\int_t^T|Y_i^n(s)|^2\, ds+\epsilon E\int_t^T|Z_i^n(s)|^2\, ds\\
& &+\epsilon
E\int_t^Tn^2\sum_{l=1}^m\left((Y_i^n(s)-Y_l^n(s)-k(i,l))^+\right)^2\, ds\\
&\le &\displaystyle C_\epsilon +C_\epsilon E\int_t^T|Y_i^n(s)|^2\,
ds+C\epsilon E\int_t^T|Z_i^n(s)|^2\, ds. \end{eqnarray*} Here,
$C_\epsilon>0$ denotes a constant which depends on $\epsilon$ and
may vary from line to line.

 Therefore,
 \be\ba{c} \displaystyle
 E|Y_i^n(t)|^2+E\int_t^T|Z_i^n(s)|^2\, ds\le C.
\ea\ee

From (\ref{nn}), we obtain (\ref{estimateK}), and  from
(\ref{finalito}), we deduce
 \be\ba{c} \displaystyle
 ||Y^n||_{S^2}
\le C. \ea \ee The proof of Lemma \ref{apriori} is now complete.
\endpf

\medskip

\section{Uniqueness}

In this section, we prove the uniqueness by a verification method.
Let $(\wt Y,\wt Z,\wt K)$ be a solution in the space
$(S^2,M^2,N^2)$ to RBSDE (\ref{RBSDEi}). We will prove that $\wt
Y$ is in fact the (vector) value for an optimal switching problem
of BSDEs. For this purpose, we introduce the following optimal
switching problem.

Let $\{\theta_j\}_{j=0}^\infty$ be an increasing sequence of
stopping times with values in $[0,T]$ and $\forall j$, $\alpha_j$
is an ${\cal F}_{\theta_j}$-measurable random variable with values
in $\Lambda$, and $\chi$ is the indicator function. We assume
moreover that for $P$-a.s. $\omega$, there exists an integer
$N(\omega)$ such that $\theta_{N}=T$.

Then we define the admissible switching strategy as follows:
\begin{equation}\label{admi}
a(s)=\alpha_0\chi_{\{\theta_0\}}(s)+\sum_{j=1}^N
\alpha_{j-1}\chi_{(\theta_{j-1},\theta_j]}(s).
\end{equation}

We denote by ${\cal A}$ the set of all these admissible switching
strategies and by ${\cal A}^i$ the subset of ${\cal A}$ consisting
of admissible switching strategies starting from the mode $i$. In
the same way, we denote by ${\cal A}_t$ the set of all the
admissible strategies starting at the time $t$ (or equivalently
$\theta_0=t$ ) and by ${\cal A}^i_t$ the subset of ${\cal A}_t$
consisting of admissible switching strategies starting at time $t$
from the mode $i$.

For any $a(\cdot)\in {\cal A}_t$, we define the associated (cost)
process $A^{a(\cdot)}$ as follows:
\begin{equation}\label{UA}
A^{a(\cdot)}(s)=\sum_{j=1}^{N-1}
k(\alpha_{j-1},\alpha_j)\chi_{[\theta_{j},T]}(s),\ s\in [t,T].
\end{equation}
Obviously, $A^{a(\cdot)}(\cdot)$ is a c\`adl\`ag process.

Now we are in position to introduce the switched BSDE: \be
\label{UV} U(s) =\xi_{a(T)}
+A^{a(\cd)}(T)-A^{a(\cd)}(s)+\int_s^T\psi(r,U(r),V(r),a(r))
dr-\int_s^TV(r)dW(r),s\in [t,T].\ee This is a (slightly)
generalized BSDE: it is equivalent to the following standard BSDE:
\be \label{UV1} {\bar U}(s) =\xi_{a(T)}
+A^{a(\cd)}(T)+\int_s^T\psi(r,\bar{U}(r)-A^{a(\cdot)}(r),\bar{V}(r),a(r))
dr-\int_s^T\bar{V}(r)dW(r),s\in [t,T]\ee via the simple change of
variable:
$$\bar{U}(s)=U(s)+A^{a(\cdot)}(s),\quad \bar{V}(s)=V(s),\quad s\in [t,T].$$

Hence, BSDE (\ref{UV}) has a solution in $S^2\times M^2$. We
denote this solution by $(U^{a(\cdot)},V^{a(\cdot)})$. Note that
$U$ is only a c\`adl\`ag process.

\medskip

The optimal switching problem with the initial mode $i\in\Lambda$
is to minimize $U^{a(\cdot)}(t)$ subject to $a(\cdot)\in {\cal
A}_t^i$.

\medskip

The assumptions  required for the uniqueness will be slightly
stronger than those needed for existence. We keep the same
assumption on $\psi$ and we assume the following for $k$.

\begin{hypothesis}\label{k'}
(i) For any $(i,j)\in \Lambda\times \Lambda$, $k(i,j)\ge 0$.

(ii) For any $(i,j,l)\in \Lambda\times\Lambda\times \Lambda$ such
that $i\not=j$ and $j\not=l$,
$$k(i,j)+k(j,l)> k(i,l).$$
\end{hypothesis}

We have the following representation for the first component of
the adapted solution to RBSDE (\ref{RBSDEi}), which immediately
implies the uniqueness of the adapted solution to RBSDE
(\ref{RBSDEi}).

\begin{theorem}\label{uniqueness} Let us suppose that the
Hypotheses \ref{psi} and \ref{k'} hold. Let us also assume that
$\xi\in L^2(\Omega,{\cal F}_T,P;R^m)$ takes values in $\bar{Q}$.
Let $(\wt Y,\wt Z,\wt K)$ be a solution in $(S^2,M^2,K^2)$ to
RBSDE (\ref{RBSDEi}). Then

(i) For any $a(\cdot)\in {\cal A}_t^i$, we have:
\begin{equation}\label{inequality}\wt Y_i(t)\le U^{a(\cdot)}(t), \quad P-a.s.
\end{equation}

(ii) Set $\theta_0^*=t$, $\alpha_0^*=i$. We define the sequence
$\{\theta_j^*,\alpha_j^*\}_{j=1}^{\infty}$ in an inductive way as
follows:
\begin{equation}\label{theta*}
\theta^*_j:=\inf\{s\ge \theta^*_{j-1}:
\tilde{Y}_{\a^*_{j-1}}(s)=\min_{l\not={\a^*_{j-1}}}\{\tilde{Y}_l(s)+k(\alpha^*_{j-1},l)\}\wedge
T,
\end{equation}
and $\alpha_j^*$ is the ${\cal F}_{\theta_j^*}$-measurable random
variable such that
$$\wt Y_{\a^*_{j-1}}(\theta_j^*)=\wt Y_{\alpha_j^*}(\theta^*_j)+k(\a_{j-1}^*,\a^*_j),$$
with $j=1,2,\cdots.$

Then, $P$-a.s. $\omega$, there exists an integer $N(\omega)$ such
that $\theta^*_{N}=T$. And the following  switching strategy:
\begin{equation}
a^*(s)=i\chi_{\{t\}}(s)+\sum_{j=1}^N
\alpha_{j-1}^*\chi_{(\theta^*_{j-1},\theta^*_j]}(s),
\end{equation}
is admissible, i.e., $a^*(\cdot)\in {\cal A}^i_t$. Moreover,
$$\wt Y_i(t)=U^{a^*(\cdot)}(t).$$

(iii) We have the following representation for $\wt
Y(t)$:
$$\wt Y_i(t)=\mathop{\mbox{\rm essinf}}_{a(\cdot)\in {\cal
A}_t^i} U^{a(\cdot)}(t), \quad i\in \Lambda,t\in [0,T].$$ RBSDE
(\ref{RBSDEi}) has a unique solution.
\end{theorem}

{\bf Proof.}  Without loss of generality, we will prove (i) and
(ii) for the case of $t=0$. Otherwise, it suffices to consider the
admissible switching strategies starting at time $t$.

(i) We define \begin{eqnarray} \wt Y^{a(\cd)}(s)&=&\sum_{i=1}^N\wt
Y_{\a_{i-1}}(s)\chi_{[\th_{i-1},
\th_i)}(s)+\xi_{a(T)}\chi_{\{T\}}(s),\label{UY}\\
\wt Z^{a(\cd)}(s)&=&\sum_{i=1}^N\wt
Z_{\a_{i-1}}(s)\chi_{[\th_{i-1}, \th_i)}(s),\label{UZ}\\
\wt K^{a(\cdot)}(s)&=&\sum_{i=1}^N \int_{\th_{i-1}\wedge
s}^{\th_i\wedge s}d\wt K_{\alpha_{i-1}}(r).\label{UK}
 \end{eqnarray}

Noting that $\wt Y^{a(\cd)}(\cd)$ is a c\`adl\`ag process with
jump $\wt Y_{\alpha_i}(\theta_i)-\wt Y_{\alpha_{i-1}}(\theta_i)$
at $\theta_i$, $i=1,\cdots,N-1$, we deduce that
\begin{eqnarray*}
& &\wt Y^{a(\cd)}(s)-\wt Y^{a(\cd)}(0)\\
&=&\sum_{i=1}^N \int_{\th_{i-1}\wedge s}^{\th_i\wedge
s}[-\psi(r,\wt Y_{\a_{i-1}}(r),\wt Z_{\a_{i-1}}(r),\a_{i-1})dr+\wt
Z_{\alpha_{i-1}}(r)dW(r)+d\wt
K_{\alpha_{i-1}}(r)]\\
& &+\sum_{i=1}^{N-1}[\wt Y_{\alpha_i}(\theta_i)-\wt
Y_{\alpha_{i-1}}(\theta_i)]\chi_{[\theta_i,T]}(s)\\
&=&\int_0^s[-\psi(r,\wt Y^{a(\cd)}(r),\wt
Z^{a(\cd)}(r),a(r))dr+\wt Z^{a(\cd)}(r)dW(r)+ d\wt
K^{a(\cdot)}(r)]\\
& &+\tilde{A}^{a(\cdot)}(s)-A^{a(\cdot)}(s),
\end{eqnarray*}
where \be\label{UAtilde}\wt A^{a(\cdot)}(s)=\sum_{i=1}^{N-1}[\wt
Y_{\alpha_i}(\theta_i)+k(\alpha_{i-1},\alpha_i)-\wt
Y_{\alpha_{i-1}}(\theta_i)]\chi_{[\theta_i,T]}(s),\ee and it is an
increasing process due to the fact that
$$\wt Y(t)\in \bar{Q},\quad \forall t\in [0,T].$$
 Consequently, we conclude that $(\wt Y^{a(\cd)},\wt
Z^{a(\cd)})$ is a solution of the following BSDE:
\begin{eqnarray}\label{bsdetilde}
& &\wt Y^{a(\cd)}(s)\nonumber\\
&=&\xi_{a(T)}+A^{a(\cdot)}(T)-A^{a(\cdot)}(s)-[(\wt
K^{a(\cdot)}(T)+\wt A^{a(\cdot)}(T))-
(\wt K^{a(\cdot)}(s)+\wt A^{a(\cdot)}(s))]\nonumber \\
& &+\int_s^T\psi(r,\wt Y^{a(\cd)}(r),\wt
Z^{a(\cd)}(r),a(r))dr-\int_s^T\wt Z^{a(\cd)}(r)dW(r), s\in [0,T].
\end{eqnarray}

Since both $\wt K^{a(\cdot)}$ and $\wt A^{a(\cdot)}$ are
increasing c\`adl\`ag processes, from the comparison theorem, we
conclude that
$$\wt Y^{a(\cdot)}(0)\le U^{a(\cdot)}(0),$$
which implies that
$$\wt Y_i(0)\le U^{a(\cdot)}(0).$$

(ii) Let us first claim that if $0\le \theta_1^*<\theta_2^*<T$,
then there exists a constant $c>0$ such that
$$|\wt Y(\theta_2^*)-\wt Y(\theta_1^*)|\ge c.$$

To prove this claim, we introduce the following subsets of
$\bar{Q}$: for $i\not=j$,
$$B_{i,j}:=\{(y_1,\cdots,y_m)^T\in R^m\ :\ y_i=y_j+k(i,j)\}\cap \bar{Q}.$$

We assert that for $i\not=j$ and $j\not=l$,  $B_{i,j}\cap
B_{j,l}=\emptyset$.

In fact, if there exists an element $(y_1,\cdots,y_m)^T \in
B_{i,j}\cap B_{j,l}$, then $$y_i=y_j+k(i,j) \mbox{  and }
y_j=y_l+k(j,l).$$ We deduce then
$$y_i=y_l+k(i,j)+k(j,l)>y_l+k(i,l),$$ which contradicts the fact
that  $(y_1,\cdots,y_m)^T\in \bar{Q}$.

Hence, the distance between $B_{i,j}$ and $B_{j,l}$ is strictly
positive,
$$dist(B_{i,j},B_{j,l})>0.$$

We set
$$c:=\min_{i\not=j,j\not=l}dist(B_{i,j},B_{j,l})>0.$$

We return  to the proof of the claim. From the definition of
$(\theta^*_1, \alpha^*_1)$ and $(\theta^*_2, \alpha^*_2)$,
$$\wt Y(\theta_1^*)\in B_{i, \a_1^*} \mbox{ and } \wt Y(\theta_2^*)\in B_{\a_1^*,
\a_2^*},$$ which implies that
$$|\wt Y(\theta_2^*)-\wt Y(\theta_1^*)|\ge c,$$
and the proof of the claim is finished.

In the same way, if
$\theta_1^*<\theta_2^*<\cdots<\theta_{j-1}^*<\theta_j^*<T$, then
$$|\wt Y(\theta_j^*)-\wt Y(\theta_{j-1}^*)|\ge c.$$

On the other hand, as $\wt Y$ satisfies (\ref{RBSDEi}), it is easy
to check that
$$E\left[\sum_{j=1}^\infty
|\wt Y(\theta_j^*)-\wt Y(\theta_{j-1}^*)|^2\right]<\infty.$$

As a  consequence, there exists $N(\omega)$ such that
$\theta_N^*=T$.

Finally, from the choice of $a^*(\cdot)$,
$$\wt K^{a^*(\cdot)}+\wt A^{a^*(\cdot)}=0.$$
We conclude from (\ref{bsdetilde}) that
$$\wt Y^{a^*(\cdot)}(0)= U^{a^*(\cdot)}(0),$$
which implies that
$$\wt Y_i(0)= U^{a^*(\cdot)}(0).$$

(iii) The representation for $\wt Y$ is a combination of both
assertions (i) and (ii). This gives  the uniqueness of the first
component of the adapted solution, and the uniqueness of the other
two components of the adapted solution follows then.
\endpf

\section{Optimal switching of functional SDEs}

In this section, we  study the optimal switching problem. In order
to ensure the existence of optimal switching strategy, we use the
weak formulation of the problem. Let $(\Omega, {\cal F}, P)$ be a
complete probability space and let $\{{\cal F}_t, t\ge 0\}$ be a
filtration satisfying the usual conditions. The process $W$ is an
$\{{\cal F}_t, t\ge 0\}$-Brownian motion on $R^d$ defined on
$(\Omega, {\cal F},P)$.

 Consider
the switched equation
\begin{equation}\label{sdef2}
X^{a(\cdot)}(t)=x_0+\int_0^t
\sigma(s,X^{a(\cdot)})[dW(s)+b(s,X^{a(\cdot)},a(s))ds], \quad t\in
[0,T]
\end{equation}
and the cost functional
\begin{equation}\label{cost2}
J(a(\cdot))=E\left[\int_0^T
l(s,X^{a(\cdot)},a(s))ds\right]+E\left[\sum_{i=1}^{N-1}
k(\alpha_{i-1},\alpha_i)\right].
\end{equation}
The switching problem is to minimize the cost $J(a(\cdot))$ with
respect to $a(\cdot)$, subject to the state equation
(\ref{sdef2}).

 In the above, $x_0$ is a fixed point in $R^d$. $\sigma$, $b$ and $l$ are defined on
$[0,T]\times C([0,T];R^d)$, $[0,T]\times C([0,T];R^d)\times
\Lambda$ and $[0,T]\times C([0,T];R^d)\times \Lambda$,
respectively, with values in $R^{d\times d}$, $R^d$ and $R$,
respectively. As in Section 3, $a(\cdot)$ is an admissible $({\cal
F}_t)_{t\ge 0}$-adapted switching strategy, ${\cal A}$ is the set
of all the admissible $\{{\cal F}_t, t\ge 0\}$-adapted switching
strategies and ${\cal A}^i$ is the subset of ${\cal A}$ consisting
of  the admissible $\{{\cal F}_t, t\ge 0\}$-adapted switching
strategies starting from the mode $i$. We assume that $\sigma$,
$b(\cdot,i)$ and $l(\cdot,i)$ are progressively measurable
functionals on $C([0,T];R^d)$ in the following sense:
\begin{definition} Let $C([0,T];R^d)$ be the space of continuous
functions $x:[0,T]\rightarrow R^d$.  For $0\le t\le T$, define
${\cal G}_t:=\sigma(x(s):0\le s\le t)$, and set ${\cal G}:={\cal
G}_T$. A progressively measurable functional on $C([0,T];R^d)$ is
a mapping $\mu:[0,T]\times C([0,T];R^d)\to H$ ( $H$ is some
Euclidean space) such that for each fixed $t\in [0,T]$, $\mu$
restricted to $[0,t]\times C([0,T];R^d)$ is ${\cal
B}([0,t])\otimes {\cal G}_t/{\cal B}(H)$-measurable.
\end{definition}

We assume that $k$ satisfies Hypothesis \ref{k'}. And we assume
also that $\sigma$, $b$ and $l$ satisfy the following hypothesis.
\begin{hypothesis}\label{sbl}
(i) $\sigma$, $b(\cdot,\cdot,i)$ and $l(\cdot,\cdot,i),i\in
\Lambda,$ are progressively measurable functionals on
$C([0,T];R^d)$.

(ii) There exists a constant $\beta
>0$ such that $\forall (t,x,x',i)\in [0,T]\times R^d\times R^d\times\Lambda$,
$$|b(t,x,i)-b(t,x',i)|+|\sigma(t,x)-\sigma(t,x')|+|l(t,x,i)-l(t,x',i)|\le\beta
||x-x'||_{C([0,t];R^d)}.$$

(iii) $\sigma$ has a bounded inverse.

(iv) $b$ is bounded.
\end{hypothesis}

Let $(Y,Z,K)$ be the unique solution in $(S^2,M^2,N^2)$ of the
following  RBSDE:
\be\label{RBSDEipb}\left\{\ba{rcl}Y_i(t)&=&\displaystyle\int_t^T\psi(s,X,Z_i(s),i)\,
ds-\int_t^TdK_i(s)-\int_t^TZ_i(s)\, dW(s),\\
Y_i(s)&\le&\displaystyle \min_{j\not=i}\{Y_j(s)+k(i,j)\}, \\
 &&\displaystyle
\int_0^T\left(Y_i(s)-\min_{j\not=i}\{Y_j(s)+k(i,j)\}\right)dK_i(s)=0,
\quad i\in \Lambda,\ea\right.\ee where $\psi$ is defined as
follows: $\forall (t,x,z,i)\in [0,T]\times C([0,T];R^d)\times
R^d\times \Lambda$,
$$\psi(t,x,z,i):=l(t,x,i)+\langle z,b(t,x,i)\rangle,$$
and $X$ is the solution to the following functional SDE:
\begin{equation}\label{sdef3}
X(t)=x_0+\int_0^t \sigma(s,X)dW(s), \quad t\in [0,T].
\end{equation}

\begin{theorem}\label{theormswitch}
Let the Hypotheses \ref{k'} and \ref{sbl} hold. Then

(i) For any $a(\cdot)\in {\cal A}^i$, we have:
\begin{equation}\label{inequality2}J(a(\cdot))\ge Y_i(0).
\end{equation}

(ii) There exists an optimal switching strategy $a^*$, and  a weak
solution $(P^*, W^*,X^*)$, such that
\begin{equation}\label{weak}X^*(t)=x_0+\int_0^t
\sigma(s,X^*)[dW^*(s)+b(s,X^*,a^*(s,X^*))],  \quad t\in [0,T],
\end{equation} and
$$J(a^*(\cdot))= Y_i(0).$$
\end{theorem}

{\bf Proof.} (i) For any $a(\cdot)\in {\cal A}^i$, we set
$$d\bar{P}:=\exp\left\{-\int_0^T b(s,X^{a(\cdot)},a(s))dW(s)-\frac{1}{2}\int_0^T
|b(s,X^{a(\cdot)},a(s))|^2ds\right\}dP.$$ Then
$\bar{W}(t)=W(t)+\int_0^t b(s,X^{a(\cdot)},a(s))ds$ is a Brownian
motion under the new probability measure $\bar{P}$. Let
$(\bar{Y},\bar{Z},\bar{K})$ be the solution of the following
RBSDE:
\be\label{RBSDEbar}\left\{\ba{rcl}\bar{Y}_i(t)&=&\displaystyle\int_t^T\psi(s,X^{a(\cdot)},\bar{Z}_i(s),i)\,
ds-\int_t^Td\bar{K}_i(s)-\int_t^T\bar{Z}_i(s)\, d\bar{W}(s),\\
\bar{Y}_i(s)&\le&\displaystyle \min_{j\not=i}\{\bar{Y}_j(s)+k(i,j)\}, \\
 &&\displaystyle
\int_0^T\left(\bar{Y}_i(s)-\min_{j\not=i}\{\bar{Y}_j(s)+k(i,j)\}\right)d\bar{K}_i(s)=0,
\quad i\in\Lambda.\ea\right.\ee Note that since $X^{a(\cdot)}$
solves (\ref{sdef2}), we have
$$X^{a(\cdot)}(t)=x_0+\int_0^t \sigma(s,X^{a(\cdot)})d\bar{W}(s), \quad t\in [0,T].$$
 By a classical argument of
Yamada-Watanabe, for RBSDE (\ref{RBSDEipb}), the pathwise
uniqueness implies  the uniqueness in the sense of probability law
(see, e.g. \cite{Delarue}, for a proof in the framework of BSDE).
Hence, we have \be\label{F2}Y_i(0)=\bar{Y}_i(0), \quad i\in
\Lambda.\ee

Recalling the cost process $A^{a(\cdot)}$ defined by (\ref{UA}),
and defining  $\bar{Y}^{a(\cdot)}$, $\bar{Z}^{a(\cdot)}$ ,
$\bar{K}^{a(\cdot)}$ and $\bar{A}^{a(\cdot)}$ in the same way as
(\ref{UY}), (\ref{UZ}), (\ref{UK}) and (\ref{UAtilde}), we deduce
in the same manner as  in Section 3 that $(\bar{Y}^{a(\cdot)},
\bar{Z}^{a(\cdot)})$ is the unique solution of the following
BSDE:\begin{eqnarray}\label{bsdebar}
& &\bar{Y}^{a(\cd)}(t)\nonumber\\
&=&A^{a(\cdot)}(T)-A^{a(\cdot)}(t)-[(\bar{K}^{a(\cdot)}(T)+\bar{A}^{a(\cdot)}(T))-
(\bar{K}^{a(\cdot)}(t)+\bar{A}^{a(\cdot)}(t))]\nonumber \\
& &+\int_t^T\psi(s,X^{a(\cdot)},
\bar{Z}^{a(\cd)}(s),a(s))ds-\int_t^T\bar{Z}^{a(\cd)}(s)d\bar{W}(s),\quad
t\in [0,T].
\end{eqnarray}
Since, both $\bar{K}$ and $\bar{A}$ are increasing, we deduce that
\begin{eqnarray}\label{F1}
& &\bar{Y}^{a(\cd)}(0)\nonumber\\
&\le&A^{a(\cdot)}(T)+\int_0^T\psi(s,X^{a(\cdot)},
\bar{Z}^{a(\cd)}(s),a(s))ds-\int_0^T\bar{Z}^{a(\cd)}(s)[dW(s)+b(s,X^{a(\cdot)},a(s))]\nonumber\\
&=&A^{a(\cdot)}(T)+\int_0^T
l(s,X^{a(\cdot)},a(s))ds-\int_0^T\bar{Z}^{a(\cd)}(s)dW(s).
\end{eqnarray}
From the definition of $\bar{Y}^{a(\cdot)}$,
\be\label{F3}\bar{Y}_i(0)=\bar{Y}^{a(\cdot)}(0).\ee From
(\ref{F2}), (\ref{F3}) and (\ref{F1}),
$$Y_i(0)\le A^{a(\cdot)}(T)+\int_0^T
l(s,X^{a(\cdot)},a(s))ds-\int_0^T\bar{Z}^{a(\cd)}(s)dW(s).$$
Taking expectation with respect to $P$, we have
$$Y_i(0)\le J(a(\cdot)).$$

(ii) Let $X$ be the solution of SDE (\ref{sdef3}) and  $(Y,Z,K)$
be the solution of RBSDE (\ref{RBSDEipb}). Then, $Y$ is adapted to
the filtration ${\cal F}_t^W={\cal F}_t^X$, due to (\ref{sdef3})
and Hypothesis \ref{sbl} (iii).

Set $\theta_0^*=0$, $\alpha_0^*=i$. We define the sequence
$\{\theta_j^*,\alpha_j^*\}_{j=1}^{\infty}$ in an inductive way as
follows:
\begin{equation}
\theta^*_j:=\inf\{s\ge \theta^*_{j-1}:
{Y}_{\alpha^*_{j-1}}(s)=\min_{l\not=\alpha^*_{j-1}}\{{Y}_l(s)+k(\alpha^*_{j-1},l)\}\wedge
T,
\end{equation}
and $\alpha_j^*$ is the ${\cal F}_{\theta_j^*}$-measurable random
variable such that
$${Y}_{\alpha^*_{j-1}}(\theta_j^*)={Y}_{\alpha_j^*}(\theta^*_j)+k(\a_{j-1}^*,\a^*_j),$$
with $j=1,2,\cdots.$ Then, $P$-a.s. $\omega$, there exists an
integer $N(\omega)$ such that $\theta^*_{N}=T$. And we define the
switching strategy $a^*$ as follows:
\begin{equation}
a^*(s):=i\chi_{\{0\}}(s)+\sum_{j=1}^N
\alpha_{j-1}^*\chi_{(\theta^*_{j-1},\theta^*_j]}(s),
\end{equation}
is admissible, i.e., $a^*(\cdot)\in {\cal A}^i$. $a^*$ is adapted
to the filtration ${\cal F}_t^W={\cal F}_t^X$, i.e.,
$a^*(t)=a^*(t,X).$ Setting $$dP^*:=\exp\left(\int_0^t
b(s,X,a^*(s,X))dW(s)-\frac{1}{2} \int_0^t
|b(s,X,a^*(s,X))|^2ds\right)dP,$$ then
$$W^*(t):=W(t)-\int_0^t b(s,X,a^*(s,X))dW(s), \quad t\in [0,T],$$
is a Brownian motion under the probability measure $P^*$, and
$(P^*,W^*,X)$ is a weak solution of (\ref{sdef2}).

Computing $Y^{a^*(\cdot)}$ as in Section 3, we deduce that
\begin{eqnarray}\label{bsde4}
& &Y^{a^*(\cd)}(0)\nonumber\\
&=&A^{a^*(\cdot)}(T)-( K^{a^*(\cdot)}(T)+ \wt A^{a^*(\cdot)}(T))\nonumber \\
& &+\int_0^T\psi(s, X, Z^{a^*(\cd)}(s),a^*(s))ds-\int_0^T
Z^{a^*(\cd)}(s)dW(s),
\end{eqnarray}
where \begin{eqnarray*} A^{a^*(\cdot)}(s)&=&\sum_{j=1}^{N-1}
k(\alpha^*_{j-1},\alpha^*_j)\chi_{[\theta^*_{j},T]}(s), \\
Y^{a^*(\cd)}(s)&=&\sum_{i=1}^N
Y_{\a^*_{i-1}}(s)\chi_{[\th^*_{i-1},
\th^*_i)}(s)+Y_{a^*(T)}(s)\chi_{\{T\}}(s),\\
 Z^{a^*(\cd)}(s)&=&\sum_{i=1}^N
Z_{\a^*_{i-1}}(s)\chi_{[\th^*_{i-1}, \th^*_i)}(s),\\
K^{a^*(\cdot)}(s)&=&\sum_{i=1}^N \int_{\th^*_{i-1}\wedge
s}^{\th_i^*\wedge s}dK_{\alpha^*_{i-1}}(r),\\
\tilde{A}^{a^*(\cdot)}(s)&=&\sum_{i=1}^{N-1}[
Y_{\alpha_i}(\theta^*_i)+k(\alpha^*_{i-1},\alpha^*_i)-
Y_{\alpha^*_{i-1}}(\theta^*_i)]\chi_{[\theta_i^*,T]}(s).
\end{eqnarray*}
From the definition of $a^*(\cdot)$, $K^{a^*(\cdot)}=0$ and
$A^{a^*(\cdot)}=0$. Hence, from (\ref{bsde4}), it follows that
\begin{eqnarray}\label{cost4}
& &Y^{a^*(\cd)}(0)\nonumber\\
&=&A^{a^*(\cdot)}(T) +\int_0^T l(s, X, a^*(s))ds-\int_0^T
Z^{a^*(\cd)}(s)dW^*(s).\end{eqnarray} Taking the expectation with
respect to $P^*$, we conclude the proof.
\endpf

\section{System of variational inequalities}

In this section, we will show that the RBSDE studied in Sections 2
and 3 allows us to give a probabilistic representation of the
solution to a system of variational inequalities. For this
purpose, we will put  RBSDE (\ref{RBSDEi}) in a Markovian
framework.

\medskip

Let $b: [0,T]\times R^d\rightarrow R^d$ and $\sigma:[0,T]\times
R^d\rightarrow R^{d\times d}$ be continuous mappings. We assume:
\begin{hypothesis}\label{bsigma}
There exists a constant $C>0$, such that for all $t\in [0,T]$,
and $(x,x')\in R^d\times R^d$,
$$|b(t,0)|+|\sigma(t,0)|\le C,$$
and
$$|b(t,x)-b(t,x')|+|\sigma(t,x)-\sigma(t,x')|\le C|x-x'|.$$
\end{hypothesis}

\medskip

For each $(t,x)\in [0,T]\times R^d$, let $\{X^{t,x}(s); t\le s\le
T\}$ be the unique $R^d$-valued solution of the SDE:
$$X^{t,x}(s)=x+\int_t^s
b(r,X^{t,x}(r))dr+\int_t^s\sigma(r,X^{t,x}(r))dW(r),\quad s\in
[t,T].$$

We suppose now that the data $\xi$ and $\psi$ of RBSDE
(\ref{RBSDEi}) take the following form:
\begin{eqnarray*}
\xi_i&=&g(X^{t,x}(T),i),\\
\psi(s,y,z,i)&=&\psi(s,X^{t,x}(s),y,z,i),
\end{eqnarray*}
where $g$ and $\psi$ are given as follows.
\begin{hypothesis}\label{gpsi}
(i) For each $i\in \Lambda$, the function $g(\cdot,i)\in C(R^d)$
and has at most polynomial growth at infinity.

(ii) For each $i\in \Lambda$,
$$\psi(\cdot,i):[0,T]\times R^d\times R\times R^d\rightarrow R$$ is
jointly continuous and there exist two constants $C>0$ and $p\ge
0$ such that
$$|\psi(t,x,0,0,i)|\le C(1+|x|^p),$$
$$|\psi(t,x,y,z,i)-\psi(t,x,y',z',i)|\le C(|y-y'|+|z-z'|),$$
for $t\in [0,T],x,z,z'\in R^d, y,y'\in R$, and $i\in \Lambda$.

(iii) $\forall x\in R^d$,
$$g(x):=(g(x,1),\cdots,g(x,m))^T\in
\bar{Q}.$$
\end{hypothesis}

For each $t\ge 0$, we denote by $\{{\cal F}_s^t,t\le s\le T\}$ the
natural filtration of the Brownian motion $\{W_s-W_t,t\le s\le
T\}$, augmented by the $P$-null sets of ${\cal F}$.

It follows from the results of Sections 2 and 3 that for each
$(t,x)$, there exists a unique triple $(Y^{t,x},Z^{t,x},K^{t,x})$
in $S^2\times M^2\times N^2$ of $\{{\cal F}_s^t, t\le s\le T\}$
progressively measurable processes, which solves the following
RBSDE: \be\label{RBSDE5}\left\{\ba{rcl}Y_i(s)&=&\displaystyle
g(X^{t,x}(T),i)+\int_s^T\psi(r,X^{t,x}(r),Y_i(r),Z_i(r),i)
dr\\
& &\displaystyle-\int_s^TdK_i(r)-\int_s^TZ_i(r)\, dW(r),\ s\in [0,T];\\
Y_i(s)&\le&\displaystyle \min_{j\not=i}\{Y_j(s)+k(i,j)\},\ s\in [0,T]; \\
 &&\displaystyle
\int_0^T\left(Y_i(s)-\min_{j\not=i}\{Y_j(s)+k(i,j)\}\right)dK_i(s)=0;
\ i\in\Lambda.\ea\right.\ee

\medskip

We now consider the related system of variational inequalities.
Roughly speaking, a solution of the system of variational
inequalities is a function $u: [0,T]\times R^d\rightarrow R^m$
which satisfies:
\begin{eqnarray}\label{vi5}
\max&\biggl\{&-\partial_t u_i(t,x)-{\cal
L}u_i(t,x)-\psi(t,x,u_i(t,x),\nabla
u_i(t,x)\sigma(t,x),i),\nonumber\\
& & u_i(t,x)-\min_{j\not=i}(u_j(t,x)+k(i,j))\biggr\}=0,
\end{eqnarray}
$(t,x,i)\in [0,T]\times R^d\times\Lambda$, with the terminal
condition
\begin{equation}\label{tc5}
u_i(T,x)=g(x,i),\quad (x,i)\in R^d\times\Lambda.
\end{equation}
Here, the second-order partial differential operator ${\cal L}$ is
given by
$${\cal L}:=\frac{1}{2} \sum_{j,l=1}^d
(\sigma\sigma^T(t,x))_{j,l}\frac{\partial^2}{\partial x_j\partial
x_l}+\sum_{j=1}^d b_j(t,x)\frac{\partial}{\partial x_j}.$$

More precisely, we shall consider solution of (\ref{vi5}) in the
viscosity sense. It will be convenient for the sequel to define
the notion of viscosity solution in the language of sub- and
super-jets (see, e.g., \cite{CIL}). Below, $S(d)$ will denote the
set of $d\times d$ symmetric nonnegative matrices.
\begin{definition}\label{jet}
Let $u\in C((0,T)\times R^d;R)$ and $(t,x)\in (0,T)\times R^d$. We
denote by ${\cal P}^{2,+} u(t,x)$ [the ``parabolic superjet'' of
$u$ at $(t,x)$] the set of triples $(p,q,X)\in R\times R^d\times
S(d)$ which are such that
$$u(s,y)\le
u(t,x)+p(s-t)+\langle q,y-x\rangle +\frac{1}{2}\langle
X(y-x),y-x)\rangle +o(|s-t|+|y-x|^2).$$ Similarly, we denote by
${\cal P}^{2,-} u(t,x)$ [the ``parabolic subjet'' of $u$ at
$(t,x)$] the set of triples $(p,q,X)\in R\times R^d\times S(d)$
which are such that
$$u(s,y)\ge
u(t,x)+p(s-t)+\langle q,y-x\rangle +\frac{1}{2}\langle
X(y-x),y-x)\rangle +o(|s-t|+|y-x|^2).$$
\end{definition}

We can now give the definition of a viscosity solution of the
system of  variational inequalities (\ref{vi5}) and (\ref{tc5}).

\begin{definition}\label{viscosity5}
 $u\in C([0,T]\times R^d;R^m)$ is called a viscosity subsolution
(resp., supersolution)  of (\ref{vi5}) and (\ref{tc5}) if
$u_i(T,x)\le \ (\mbox{resp.} \ge)\ g(x,i), (x,i)\in
R^d\times\Lambda$,  and at any point $(t,x,i)\in (0,T)\times
R^d\times \Lambda$, for any $(p,q,X)\in {\cal P}^{2,+} u_i(t,x)$
(resp., ${\cal P}^{2,-} u_i(t,x)$) ,
\begin{eqnarray*}\max&\biggl\{&-p-\frac{1}{2}Tr(\sigma\sigma^T(t,x)X)-\langle b(t,x),q\rangle -\psi(t,x,u_i(t,x),q\sigma(t,x),i),\\
& & u_i(t,x)-\min_{j\not=i}(u_j(t,x)+k(i,j))\biggr\}\le  \
(\mbox{resp.,} \ge)\  0. \end{eqnarray*}

$u\in C([0,T]\times R^d;R^m)$ is called a viscosity solution of
(\ref{vi5}) and (\ref{tc5}) if it is both a subsolution and a
supersolution of  (\ref{vi5}) and (\ref{tc5}).
\end{definition}

We now define
\begin{equation}\label{u5}
u_i(t,x):=Y_i^{t,x}(t),\quad (t,x,i)\in [0,T]\times R^d\times
\Lambda; \quad u:=(u_1,u_2,\cdots,u_m)^T.
\end{equation}
Note that $u$ is  deterministic.

\begin{lemma} Let the Hypotheses \ref{k'},  \ref{bsigma} and \ref{gpsi}
hold. For each $i\in \Lambda$, $u_i\in C([0,T]\times R^d;R)$.
\end{lemma}

{\bf Proof:} For $(t,x)\in [0,T]\times R^d$, $i\in\Lambda$ and
$a(\cdot)\in {\cal A}_t^i$, let
$(U^{a(\cdot)}_{t,x},V^{a(\cdot)}_{t,x})$ be the unique solution
of the following switched BSDE:
\begin{eqnarray} \label{UV5} U(s) &=&g(X^{t,x}(T),a(T))
+A^{a(\cd)}(T)-A^{a(\cd)}(s)\nonumber\\
& &+\int_s^T\psi(r,X^{t,x}(r),U(r),V(r),a(r)) dr\nonumber\\
& &-\int_s^TV(r)\, dW(r), \quad s\in [0,T].\end{eqnarray}

From Theorem \ref{uniqueness}, we have
$$u_i(t,x)=\mathop{\inf}_{a(\cdot)\in {\cal A}_t^i} U_{t,x}^{{a(\cdot)}}(t),\quad (t,x,i)\in [0,T]\times R^d\times\Lambda.$$

By some classical stability arguments, we obtain the continuity of
$u_i$.
\endpf

\begin{theorem} Let the Hypotheses \ref{k'}, \ref{bsigma} and \ref{gpsi}
be true. The  function $u$ given by (\ref{u5}) is the viscosity
solution of the system of variational inequalities (\ref{vi5}) and
(\ref{tc5}).
\end{theorem}
{\bf Proof:} We are going to approximate RBSDE (\ref{RBSDE5}) by
penalization, which was studied in Section 2. For each $(t,x)\in
[0,T]\times R^d$, let $\{^nY^{t,x}(s),^nZ^{t,x}(s), t\le s\le T\}$
denote the solution of the penalized BSDE:
\begin{equation}\label{penalized5}
 Y_i(s)= g(X^{t,x}(T),i)
+\int_s^T\psi^n(r,X^{t,x}(r),Y(r),Z_i(r),i)\, dr-\int_s^T Z_i(r)\,
dW(r),
\end{equation}
where for $(t,x,y,z,i)\in [0,T]\times R^d\times R^m\times
R^d\times\Lambda$,
$$\psi^n(t,x,y,z,i):=\psi(t,x,y_i,z,i)-n\sum_{j\not=
i}(y_i-y_j-k(i,j))^+.$$
 It is known from \cite{PP2} and \cite{BBP} that
$$u^n(t,x):=^nY^{t,x}(t), \quad (t,x)\in [0,T]\times R^d,$$
is the viscosity solution to the following system of parabolic
PDEs:
\begin{eqnarray}\label{pde5}
-\partial_t u^n_i(t,x)-{\cal L}u^n_i(t,x)&
-&\psi^n(t,x,u^n(t,x),\nabla
u^n_i(t,x)\sigma(t,x),i)=0,\nonumber\\
& &u_i^n(T,x)=g(x,i),\quad (t,x,i)\in [0,T]\times R^d\times
\Lambda.
\end{eqnarray}

However, from the results of Section 2, for each $(t,x,i)\in
[0,T]\times R^d\times\Lambda$,
$$u_i^n(t,x)\downarrow u_i(t,x), \mbox { as } n\rightarrow \infty.$$
Since $u^n$ and $u$ are continuous, it follows from Dini's theorem
that the above convergence is uniform on compacts.

We now show that $u$ is a subsolution of (\ref{vi5}) and
(\ref{tc5}). Let $(t,x,i)$ be a point in $[0,T]\times R^d\times
\Lambda$. Since $u$ is defined by (\ref{u5}),
$$u_i(t,x)\le \min_{j\not=i} (u_j(t,x)+k(i,j)),$$
and
$$u_i(T,x)=g(x,i).$$

Let  $(t,x,i)\in (0,T)\times R^d\times\Lambda$ and $(p,q,X)\in
{\cal P}^{2,+} u_i(t,x)$. From Lemma 6.1 in \cite{CIL}, there
exist sequences
\begin{eqnarray*}
n_l&\rightarrow &+\infty,\\
(t_l,x_l)&\rightarrow& (t,x),\\
(p_l,q_l,X_l)&\in & {\cal P}^{2,+} u_i^{n_l}(t,x), \end{eqnarray*}
such that
$$(p_l,q_l,X_l)\rightarrow (p,q,X).$$

On the other hand, for any $l$, from (\ref{pde5}),
$$-p_l-\frac{1}{2}Tr(\sigma\sigma^T(t,x)X_l)-\langle b(t,x),q_l\rangle-\psi(t,x,u_i^{n_l}(t,x),q_l\sigma(t,x),i)\le
0.
$$
Hence, taking the limit as $j\rightarrow \infty$ in the last
inequality yields:
$$-p-\frac{1}{2}Tr(\sigma\sigma^T(t,x)X)-\langle b(t,x),q\rangle-\psi(t,x,u_i(t,x),q\sigma(t,x),i)\le 0.$$
We have proved that $u$ is a subsolution of (\ref{vi5}).

We conclude by showing that $u$ is a supersolution of (\ref{vi5})
and (\ref{tc5}). Let $(t,x,i)\in (0,T)\times R^d\times\Lambda$ be
a point at which $u_i(t,x)<\min_{j\not=i} (u_j(t,x)+l(i,j))$, and
let $(p,q,X)\in {\cal P}^{2,-} u_i(t,x)$. Again from  Lemma 6.1 in
\cite{CIL}, there exist sequences
\begin{eqnarray*}
n_l&\rightarrow &+\infty,\\
(t_l,x_l)&\rightarrow& (t,x),\\
(p_l,q_l,X_l)&\in & {\cal P}^{2,-} u_i^{n_l}(t,x), \end{eqnarray*}
such that
$$(p_l,q_l,X_l)\rightarrow (p,q,X).$$

On the other hand, for any $l$, from (\ref{pde5}),
$$-p_l-\frac{1}{2}Tr(\sigma\sigma^T(t,x)X_l)-\langle b(t,x),q_l\rangle-\psi^{n_l}(t,x,u^{n_l}(t,x),q_l\sigma(t,x),i)\ge
0.
$$
From the assumption that $u_i(t,x)<\min_{j\not=i}
(u_j(t,x)+l(i,j))$ and the uniform convergence of $u^n$, it
follows that for $j$ large enough,
$u_i^{n_l}(t_l,x_l)<\min_{j\not=i} (u_j^{n_l}(t_l,x_l)+l(i,j))$;
hence, taking the limit as $j\rightarrow \infty$ in the last
inequality yields:
$$-p-\frac{1}{2}Tr(\sigma\sigma^T(t,x)X)-\langle b(t,x),q\rangle-\psi(t,x,u_i(t,x),q\sigma(t,x),i)\ge 0.$$
We have proved that $u$ is a supersolution of (\ref{vi5}) and
(\ref{tc5}).

Hence, $u$ is a viscosity solution of (\ref{vi5}) and (\ref{tc5}).
\endpf

\begin{remark}
Uniqueness of viscosity solution to (\ref{vi5}) and (\ref{tc5})
follows from classical arguments. See, e.g. \cite{CIL} and
\cite{TY}.
\end{remark}

\end{document}